# Exact Time Domain Solutions of 1-D Transient Dynamic Piezoelectric Problems with Nonlinear Damper Boundary Conditions

**Naum M. Khutoryansky, Vladimir Genis**

Department of Engineering Technology, Drexel University, Philadelphia, PA, USA
Email: khutoryn@drexel.edu, genisv@drexel.edu

## Abstract

Novel exact solutions of one-dimensional transient dynamic piezoelectric problems for thickness polarized layers and disks, or length polarized rods, are obtained. The solutions are derived using a time-domain Green's function method that leads to an exact analytical recursive procedure which is applicable for a wide variety of boundary conditions including nonlinear cases. A nonlinear damper boundary condition is considered in more detail. The corresponding nonlinear relationship between stresses and velocities at a current time moment is used in the recursive procedure. In addition to the exact recursive procedure that is effective for calculations, some new practically important explicit exact solutions are presented. Several examples of the time behavior of the output electric potential difference are given to illustrate the effectiveness of the proposed exact approach.

## Keywords

Piezoelectric Layer, Transient Dynamic Problems, Time Domain Solutions, Green's Function Method, Nonlinear Boundary Conditions, Nonlinear Damper, Output Voltage

## 1. Introduction

Piezoelectric materials and devices have been widely used in many technical applications. Nowadays, the coupling between electrical and mechanical behaviors is used in different devices based both on the so-called "direct piezoelectric effect" or the "converse piezoelectric effect" [1] [2] [3] [4].

Some newer relevant applications include (among others) the high voltage generation from transient dynamic impact processes in vehicles [5].

Analysis of operating electrical and mechanical parameters of such processes can be done by using various analytical and numerical methods. Although analytical approaches are limited to rather simple geometries and other restrictions (homogeneous or piecewise-homogeneous bodies, linear governing equations, etc.), they often provide exact solutions.

Analytical methods have been successfully used for many transient dynamic one-dimensional piezoelectric problems [6]-[14]. Among analytical methods for transient dynamic piezoelectric problems, the Laplace transform methods play a very significant role. They solve boundary value problems in the frequency-domain, possibly for complex frequencies, using transformed boundary conditions for a piezoelectric body. After obtaining such solutions, the transformation back to the time-domain employs special methods for the inversion of Laplace transforms. However, using the Laplace transform methods is not instrumental even for one-dimensional problems if nonlinear boundary conditions are considered. Time-domain numerical methods (e.g., finite element or finite difference methods) that can be used under such conditions usually lack precision associated with the use of analytical methods. Therefore, development of time-domain analytical or semi-analytical methods combining advantages of analytical and numerical methods can be of interest for such problems.

In this paper, a time-domain Green's function method is implemented for solution of one-dimensional transient dynamic piezoelectric problems for thickness polarized disks or length polarized rods. This method stems from a time-domain representation formulas approach for transient dynamic piezoe-lectric problems described in [15]. For one-dimensional problems with a variety of boundary conditions including nonlinear ones, this method produces exact solutions which are shown below. Such solutions can be used both for analyses of longitudinal mode, piezoelectric devices and as benchmark solutions for numerical methods of piezoelectricity.

## 2. Representation Formulas

Consider a transversely isotropic homogeneous piezoelectric material (piezoelectric element) with the $x_3$-axis as the poling direction and the $x_1 - x_2$ plane as the isotropic plane. Let this piezoelectric material occupy a disk (or a cylinder) $\Omega$ bounded in $x_3$-direction by planes $x_3 = 0$ and $x_3 = h$ where $h > 0$ is the thickness of the disk (or the length of the cylinder). Consider a uniaxial strain state or a stress stress state in $x_3$ direction when there is only one non-zero component of strain $\gamma_{33}$ or stress $\sigma_{33}$ the others being zero. We assume that the non-zero stress and strain components, and also the displacement $u_3$ and electric displacement $D_3$ in the $x_3$-direction, and the electric potential $\phi$, depend only on $x_3$ and $t$ which is usually the case for a longitudinal mode piezoelectric element [16]:

$$\sigma_{33} = \sigma_{33}(x_3, t), \quad D_3 = D_3(x_3, t),$$
$$u_3 = u_3(x_3, t), \quad \phi = \phi(x_3, t) \text{ in } \Omega \tag{1}$$

Under conditions (1), we can use the following one-dimensional constitutive equations (both for the uniaxial strain state and for the uniaxial stress state) that relate the mechanical and electrical fields in (1):

$$\sigma_{33} = Cu_{3,3} + e\phi_{,3}; \quad D_3 = eu_{3,3} - \epsilon\phi_{,3} \tag{2}$$

where coefficients are different for the uniaxial strain and uniaxial stress cases.

Then the corresponding equations of motions can be written as

$$\rho\ddot{u}_3 - Cu_{3,33} - e\phi_{,33} = b_3;$$
$$-D_{3,3} = \epsilon\phi_{,33} - eu_{3,33} = -q \quad \text{in } \Omega \tag{3}$$

where $b_3 = b_3(x_3, t)$ and $q = q(x_3, t)$ denote the body force in $x_3$-direction and electric charge.

To simplify further notations we will denote $x_3$ and derivatives with respect to $x_3$ by $x$ and the prime, respectively, and will skip subindex 3 for the elastic displacement, electric displacement and body force components presented in (3). Then system (3) becomes

$$\rho\ddot{u} - Cu'' - e\phi'' = b;$$
$$-D' = \epsilon\phi'' - eu'' = -q \quad \text{in } \Omega. \tag{4}$$

The Green's functions for vector $\{u, \phi\}$ can be obtained using concentrated impulses instead of $b$ or $q$ in (4) when $\Omega$ is substituted by the infinite media.

Since $\phi''$ can be expressed through $u''$ due to the second equation in (4), then the first equation in (4) can be presented as the one-dimensional wave equation for displacement $u$:

$$\rho\ddot{u} - C^D u'' = b - \frac{e}{\epsilon}q \tag{5}$$

where

$$C^D = C + \frac{e^2}{\epsilon}$$

is the Young's modulus measured at constant $D$.

The wave speed corresponding to Equation (5) is denoted below by

$$c = \sqrt{\frac{C^D}{\rho}}$$

The Green's function for $u$ corresponding to load $\{b = \delta(x)\delta(t), q = 0\}$ is the well-known Green's function for the one-dimensional wave Equation (5):

$$U(x, t) = \frac{1}{2\rho c} H(t - |x|/c) \tag{6}$$

where $H(t)$ is the Heaviside step function (right-continuous), i.e. $H(t) = 0$ for $t < 0$ and $H(t) = 1$ for $t \geq 0$.

The corresponding Green's function for $\phi$ is calculated using the second equation in (4):

$$U_\phi(x, t) = \frac{e}{\epsilon} U(x, t) \tag{7}$$

Based on (6) and (7), the representation formula for the displacement vector in 3-D case described in [15] reduces to the following expression for the displacement component $u_3(x,t) = u(x,t)$:

$$u(x,t) = \frac{1}{2}\Big[u(0,t-x/c) + u\big(h,t-(h-x)/c\big)\Big]$$
$$-\frac{1}{2\rho c}\int_{-\infty}^{t-x/c}\Big[\sigma(0,\tau) + \frac{e}{\epsilon}D(0,\tau)\Big]\mathrm{d}\tau$$
$$+\frac{1}{2\rho c}\int_{-\infty}^{t-(h-x)/c}\Big[\sigma(h,\tau) + \frac{e}{\epsilon}D(h,\tau)\Big]\mathrm{d}\tau \qquad (8)$$
$$+\frac{1}{2\rho c}\int_0^h\int_{-\infty}^{t-|\xi-x|/c}\Big[b(\xi,\tau) - \frac{e}{\epsilon}q(\xi,\tau)\Big]\mathrm{d}\tau\mathrm{d}\xi \quad \text{in } \Omega,$$

where $\sigma_{33}$ is denoted by $\sigma$ and it is taken into account that the outward normals to the lower and upper boundaries of the layer $0 \le x \le h$ have opposite directions.

In many practical applications, the electric volume charges are absent. Therefore, we consider henceforth only the case when $q = 0$. Then the terms related to $D$ in the above expression can be simplified since, based on Equation (4) in this case, $D(x,t)$ is spatially uniform:

$$D(x,t) = D(t). \qquad (9)$$

Due to the property (9) the representation formula (8) can be rewritten as

$$u(x,t) = \frac{1}{2}\Big[u(0,t-x/c) + u\big(h,t-(h-x)/c\big)\Big]$$
$$-\frac{1}{2\rho c}\int_{-\infty}^{t-x/c}\Big[\sigma(0,\tau) + \frac{e}{\epsilon}D(\tau)\Big]\mathrm{d}\tau$$
$$+\frac{1}{2\rho c}\int_{-\infty}^{t-(h-x)/c}\Big[\sigma(h,\tau) + \frac{e}{\epsilon}D(\tau)\Big]\mathrm{d}\tau \qquad (10)$$
$$+\frac{1}{2\rho c}\int_0^h\int_{-\infty}^{t-|\xi-x|/c}b(\xi,\tau)\mathrm{d}\tau\mathrm{d}\xi \quad \text{in } \Omega.$$

To obtain a representation formula for $\phi(x,t)$, let us consider an auxiliary function

$$\psi(x,t) = \phi(x,t) - \frac{e}{\epsilon}u(x,t) \qquad (11)$$

that has the following connection to the electric displacement:

$$D = -\epsilon\psi'.$$

According to (3), $\psi'' = -D'/\epsilon = 0$. Then, using the corresponding Green's function $|x|/2$ and Equation (11), we get a representation formula for $\psi(x,t)$ involving only boundary value of function $\psi(x,t)$ and a spatially uniform electric displacement:

$$\psi(x,t) = \frac{1}{2}\Big[\psi(0,t) + \psi(h,t)\Big] + \frac{h-2x}{2\epsilon}D(t) \quad \text{in } \Omega. \qquad (12)$$

Formulas (12) and (11) lead to the following expression for $\phi(x,t)$:

$$\phi(x,t) = \frac{1}{2}\Big[\phi(0,t) + \phi(h,t)\Big] - \frac{e}{2\epsilon}\Big[u(0,t) + u(h,t)\Big]$$
$$+\frac{e}{\epsilon}u(x,t) + \frac{h-2x}{2\epsilon}D(t) \quad \text{in } \Omega. \qquad (13)$$

After $u(x,t)$ is calculated, $\phi(x,t)$ can be determined using this calculated value, $D(t)$ and boundary values of $\phi(x,t)$.

The representation formula (10) allows us to get representation formulas for the velocity $v(x,t) = \dot{u}(x,t)$ and stress $\sigma(x,t)$. Differentiating (10) with respect to time provides the following representation formula for the velocity:

$$
\begin{aligned}
v(x,t) = &\frac{1}{2}\Big[ v(0, t-x/c) + v(h, t-(h-x)/c) \Big] \\
&- \frac{1}{2\rho c}\left[ \sigma(0, t-x/c) + \frac{e}{\epsilon} D(t-x/c) \right] \\
&+ \frac{1}{2\rho c}\left[ \sigma(h, t-(h-x)/c) + \frac{e}{\epsilon} D(t-(h-x)/c) \right] \\
&+ \frac{1}{2\rho c}\int_0^h b\big(\xi, t-|\xi-x|/c\big)\,\mathrm{d}\xi \quad \text{in } \Omega.
\end{aligned}
\tag{14}
$$

To get a representation formula for the stress we need to use the first contitutive equation from (2) (in the new notations introduced after equations (3)) and expression (13) which gives the following expression for the stress:

$$
\begin{aligned}
\sigma(x,t) &= Cu'(x,t) + e\left[ \frac{e}{\epsilon} u'(x,t) - \frac{1}{\epsilon} D(t) \right] \\
&= C^D u'(x,t) - \frac{e}{\epsilon} D(t).
\end{aligned}
\tag{15}
$$

After differentiating (10) with respect to $x$ and substituting the result into (15) we get

$$
\begin{aligned}
\sigma(x,t) = &\frac{\rho c}{2}\Big[ v(0, t-x/c) + v(h, t-(h-x)/c) \Big] \\
&+ \frac{1}{2}\Big[ \sigma(0, t-x/c) + \sigma(h, t-(h-x)/c) \Big] \\
&+ \frac{e}{2\epsilon}\Big[ D(t-x/c) + D(t-(h-x)/c) - 2D(t) \Big] \\
&+ \frac{1}{2}\int_0^h b\big(\xi, t-|\xi-x|/c\big)\,\mathrm{sgn}(x-\xi)\,\mathrm{d}\xi \quad \text{in } \Omega.
\end{aligned}
\tag{16}
$$

A representation formula for $D(x,t)$ is not needed under assumption that $q = 0$ since the electric displacement is uniform in space in this case and determined solely by the electric boundary conditions.

## 3. Boundary Equations

The velocity representation formula (14) generates two boundary equations when $x$ tends to the upper and lower boundaries of the piezoelectric element, that is, when $x$ tends to $h$ or 0:

$$
\begin{aligned}
v(h,t) = &v(0, t-\theta) + \frac{1}{\rho c}\Big[ \sigma(h,t) - \sigma(0, t-\theta) \Big] \\
&+ \frac{e}{\rho c \epsilon}\Big[ D(t) - D(t-\theta) \Big] + \frac{1}{\rho c}\int_0^h b\big(\xi, t-\theta+\xi/c\big)\,\mathrm{d}\xi,
\end{aligned}
\tag{17}
$$

$$v(0,t) = v(h,t-\theta) + \frac{1}{\rho c}\big[\sigma(h,t-\theta) - \sigma(0,t)\big]$$
$$+ \frac{e}{\rho c \epsilon}\big[D(t-\theta) - D(t)\big] + \frac{1}{\rho c}\int_0^h b(\xi, t-\xi/c)\mathrm{d}\xi \tag{18}$$

where $\theta$ denotes the time taken by the elastic wave to travel the thickness of the piezoelectric layer:

$$\theta = \frac{h}{c}.$$

Similarly, the stress representation formula (16) generates the following boundary equations:

$$\sigma(h,t) = \sigma(0,t-\theta) + \rho c\big[v(h,t) - v(0,t-\theta)\big]$$
$$+ \frac{e}{\epsilon}\big[D(t-\theta) - D(t)\big] - \int_0^h b(\xi, t-\theta+\xi/c)\mathrm{d}\xi, \tag{19}$$

$$\sigma(0,t) = \sigma(h,t-\theta) + \rho c\big[v(h,t-\theta) - v(0,t)\big]$$
$$+ \frac{e}{\epsilon}\big[D(t-\theta) - D(t)\big] + \int_0^h b(\xi, t-\xi/c)\mathrm{d}\xi. \tag{20}$$

It is easy to verify that Equations (17) and (19), though presented in different forms, are equivalent. The same is true for the pair of Equations (18) and (20). Therefore, we shall use the equations in these pairs interchangeably.

We will not work with boundary equations that can be obtained directly from the displacement representation formula (10), since it is computationally more effective to determine at first unknown boundary values of the velocity $v(x,t)$, and then calculate unknown boundary values of the displacement $u(x,t)$ by integrating the boundary velocity over time (using also an initial condition for $u(x,t)$).

We also need to consider boundary values of the expression (13) for the electric potential. It is important to emphasize that two equations obtained from (13) when $x$ tends to $h$ or to 0 are equivalent and, therefore, they are presented below as one equation:

$$\phi(h,t) - \phi(0,t) = \frac{e}{\epsilon}\big[u(h,t) - u(0,t)\big] - \frac{h}{\epsilon}D(t) \tag{21}$$

The boundary equations presented above will be used in the next section to create an exact time domain solution procedure in the case when nonlinear damper boundary conditions are sprecified.

## 4. Nonlinear Damper Boundary Conditions and Exact Solutions

Suppose that the lower end face of the piezoelectric element is fixed to a nonlinear damper. Let $F$ be a damping force acting on the lower end face which is defined by the following nonlinear relationship [17]:

$$F = -k_\alpha |v(0,t)|^\alpha \operatorname{sgn}(v(0,t)). \tag{22}$$

where $k_\alpha > 0$ is the damping constant, $\alpha > 0$ is the damping exponent, and $\operatorname{sgn}(.)$ is the signum function defined for all real numbers (including 0 where its value is also 0). If $v(0,t) \neq 0$, the direction of $F$ is opposite to $v(0,t)$. The exponent $\alpha$ has a value 1 for a linear damper, but may vary in practice in the interval $(0,2]$ [17] creating a set of possible boundary conditions at $x = 0$. We assume that the force $F$ is uniformly distributed over the lower end face. Then (22) transforms into the following nonlinear (in general) boundary condition at the lower end face:

$$\sigma(0,t) = \frac{k_\alpha}{A} |v(0,t)|^\alpha \operatorname{sgn}(v(0,t)) \tag{23}$$

where $A$ is the lower end face area.

Consider additional assumptions that will be used to get exact solutions for the damper boundary conditions based on the results of the previous section. We suppose that the values of $u, \sigma, b, \phi, D$ are defined for $-\infty < t < \infty$. In addition, let us assume henceforth that

$$u(x,t) = 0, \phi(x,t) = 0 \quad \text{if} \quad 0 < x < h, t < 0 \tag{24}$$

which means, based on (2) and (3), that $\sigma, D$ and $b$ are also zero inside the piezoelectric body at negative times. The next additional assumption is that

$$b(x,t) = 0 \tag{25}$$

inside the piezoelectric body at any time in the sense of generalized functions. This also includes the assumption that the initial conditions for the elastic displacement $u(x,t)$ are zero, as discussed in [15]. These assumptions will simplify using boundary Equations (17)-(20) for particular problems considered below.

Regarding the design of the piezoelectric element, we assume that it is a cylinder (or a rod) with two coated electrodes at $x = 0$ and $x = h$. The electrodes are considered to be of negligible thickness (from the mechanical point of view) and their deformation is neglected. The output voltage, which is defined as the electric potential difference between the lower and upper electrodes $\Delta\phi = \phi(0,t) - \phi(h,t)$, is of primary interest below.

The electric boundary condition at $x = 0$ corresponds to the grounded electrode:

$$\phi(0,t) = 0 \quad \text{if} \quad t \geq 0 \tag{26}$$

At the upper end face, the following mechanical boundary condition is used:

$$\sigma(h,t) = p(t) \quad \text{if} \quad t \geq 0 \tag{27}$$

where $p(t)$ is an applied normal stress load which is assumed to be known and negative.

The electric boundary condition at the upper end face $x = h$ is formulated as follows:

$$D(h,t) = 0 \quad \text{if} \quad t \geq 0. \tag{28}$$

So, based on (9), $D(t) = 0$.

Using the above assumptions the representation formulas (14) and (16) for the velocity and stress take the following simplified forms:

$$v(x,t) = \frac{1}{2}\Big[v(0,t-x/c)+v\big(h,t-(h-x)/c\big)\Big]$$
$$+\frac{1}{2\rho c}\Big[p\big(t-(h-x)/c\big)-\sigma(0,t-x/c)\Big] \quad \text{in } \Omega, \tag{29}$$

$$\sigma(x,t) = \frac{\rho c}{2}\Big[v(0,t-x/c)+v\big(h,t-(h-x)/c\big)\Big]$$
$$+\frac{1}{2}\Big[\sigma(0,t-x/c)+p\big(t-(h-x)/c\big)\Big] \quad \text{in } \Omega \tag{30}$$

where all the time dependent functions are equal to zero for negative times.

In the representation formulas (29) and (30), there are three unknown boundary functions $v(0,t), \sigma(0,t)$ and $v(h,t)$ first two of which are related by Equation (23). Two additional equations needed for determination of these three functions will be derived below based on (17) and (18).

After the the velocity $v(x,t)$ is determined for any particular $x$ over time, the corresponding displacement $u(x,t)$ can be obtained (due to the zero initial conditions) as

$$u(x,t) = \int_0^t v(x,\tau)\mathrm{d}\tau. \tag{31}$$

Boundary values of the displacement provide (according to (21) and (26)) the electric potential value at $x = h$:

$$\phi(h,t) = \frac{e}{\epsilon}\Big[u(h,t)-u(0,t)\Big] \tag{32}$$

## 4.1. An Exact Recursive Procedure

The solution of the above problem will be obtained by using an exact recursive procedure based on the following equations obtained from (17) and (18) under the boundary conditions (23) (26) (27) (28):

$$v(h,t) = 2v(0,t-\theta)-v(h,t-2\theta)+\frac{1}{\rho c}\Big[p(t)-p(t-2\theta)\Big], \tag{33}$$

$$v(0,t)+\frac{k_\alpha}{\rho cA}\big|v(0,t)\big|^\alpha \operatorname{sgn}\big(v(0,t)\big) = v(h,t-\theta)+\frac{1}{\rho c}p(t-\theta). \tag{34}$$

There are two unknowns $v(h,t)$ and $v(0,t)$ at each time moment $t$ in these equations. The right-hand sides of the equations are known at each time point since they contain either $p(t)$ or time-dalayed function values at $t-\theta$ that had to be determined at a previous step of the recursive process.

In order to simplify deriving next results, we need to introduce some additional notations:

$$\gamma = \frac{k_\alpha}{\rho cA}, \quad \xi = v(0,t), \quad r = v(h,t-\theta)+\frac{1}{\rho c}p(t-\theta). \tag{35}$$

Then, Equation (34) reads as

$$\xi + \gamma\big|\xi\big|^\alpha \operatorname{sgn}(\xi) = r. \tag{36}$$

Let the left-hand side of Equation (36) be denoted by $f(\xi)$. Since $\alpha > 0$, $f(\xi)$ is a continuous strictly monotonically increasing function on $(-\infty, \infty)$ ranging from $-\infty$ to $\infty$. Therefore, for any real $r$, there exists one and only one solution of Equation (36) in $(-\infty, \infty)$.

Denote by $Q_\alpha$ the operator that transforms $r$ into this solution of equation (36). Thus, $Q_\alpha$ is the left inverse operator of the nonlinear operator acting on $\xi$ in the left-hand side of Equation (36). If $\alpha = 2, 1, 1/2$ or $1/3$, the corresponding expressions of $Q_\alpha r$ are very simple for computations:

$$\begin{cases} Q_2 r = \dfrac{1}{2\gamma}\left(-1 + \sqrt{1 + 4|r|\gamma}\right)\operatorname{sgn}(r), \\[2mm] Q_1 r = \dfrac{r}{1+\gamma}, \; Q_{1/2} r = \dfrac{1}{4}\left(-\gamma + \sqrt{\gamma^2 + 4|r|}\right)^2 \operatorname{sgn}(r), \\[2mm] Q_{1/3} r = -\gamma\left[\dfrac{1}{6}\left(108r + 12\sqrt{12\gamma^3 + 81r^2}\right)^{1/3} - \dfrac{2\gamma}{\left(108r + 12\sqrt{12\gamma^3 + 81r^2}\right)^{1/3}}\right] + r. \end{cases} \tag{37}$$

The calculation of $Q_\alpha r$ for other values of $\alpha$ can effectively be implemented using a symbolic computation software like Maple [18].

With help of the inverse operator $Q_\alpha$ Equation (34) can be rewritten in the following explicit form for calculating $v(0,t)$:

$$v(0,t) = Q_\alpha\left[v(h,t-\theta) + \frac{1}{\rho c}p(t-\theta)\right]. \tag{38}$$

Equation (38) combined with (33) creates the recursive procedure that can be used directly for calculations or can lead to building explicit exact solution for vector $\{v(0,t), v(h,t)\}$ step by step over consecutive time intervals $j\theta \le t < (j+1)\theta$ $(j = 0, 1, 2, \cdots)$. In doing so, it is helpful to substitute $v(0, t-\theta)$ in (33) by its expression obtained from (38) which provides the following recursive equation for $v(h,t)$:

$$v(h,t) = 2Q_\alpha\left[v(h,t-2\theta) + \frac{1}{\rho c}p(t-2\theta)\right] \\ - \left[v(h,t-2\theta) + \frac{1}{\rho c}p(t-2\theta)\right] + \frac{1}{\rho c}p(t),$$

or, using the identity operator $I$ (that leaves unchanged the element on which it operates),

$$v(h,t) = (2Q_\alpha - I)\left[v(h,t-2\theta) + \frac{1}{\rho c}p(t-2\theta)\right] + \frac{1}{\rho c}p(t). \tag{39}$$

### 4.2. Explicit Exact Solutions

Now we derive some explicit exact solutions for $v(h,t)$ and $v(0,t)$ corresponding to three practically important ranges of the duration $t_1$ of the stress load at $x = h$. Our goal is to present the boundary velocities directly through the transient stress load at $x = h$ that generates the dynamic process in the piezoelectric body.

### 4.2.1. Case 1: $t_1 < 2\theta$

So, $p(t) = 0$ if $t \notin [0, 2\theta)$. Using the recursive Equation (39) under this condition for consecutive intervals $[2k\theta, 2(k+1)\theta)$, $k = 0, 1, 2, \cdots$, we finally obtain the following explicit expression for $v(h, t)$:

$$v(h,t) = \begin{cases} \dfrac{1}{\rho c} p(t) & \text{if } 0 \leq t < 2\theta, \\[2mm] (2Q_\alpha - I)^k \left[ \dfrac{2}{\rho c} p(t - 2k\theta) \right] \\[2mm] \quad \text{if } 2k\theta \leq t < 2(k+1)\theta, \ k = 1, 2, \cdots. \end{cases} \tag{40}$$

Substituting (40) into (38) we get the corresponding explicit expression for $v(0, t)$:

$$v(0,t) = \begin{cases} 0 & \text{if } 0 \leq t < \theta, \\[2mm] Q_\alpha (2Q_\alpha - I)^{k-1} \left[ \dfrac{2}{\rho c} p(t - 2k\theta) \right] \\[2mm] \quad \text{if } (2k-1)\theta \leq t < (2k+1)\theta, \ k = 1, 2, \cdots. \end{cases} \tag{41}$$

### 4.2.2. Case 2: $t_1 < 4\theta$

In this case, $p(t) = 0$ if $t \notin [0, 4\theta)$. Acting similarly to section 4 we derive the following explicit exact solutions:

$$v(h,t) = \begin{cases} \dfrac{1}{\rho c} p(t) & \text{if } 0 \leq t < 2\theta, \\[2mm] (2Q_\alpha - I)^{k-1} \left[ \dfrac{1}{\rho c} p(t - 2(k-1)\theta) \right. \\[2mm] \left. + (2Q_\alpha - I) \left( \dfrac{2}{\rho c} p(t - 2k\theta) \right) \right] \text{ if } 2k\theta \leq t < 2(k+1)\theta, k = 1, 2, \cdots; \end{cases} \tag{42}$$

$$v(0,t) = \begin{cases} 0 & \text{if } 0 \leq t < \theta, \\[2mm] Q_\alpha \left[ \dfrac{2}{\rho c} p(t - \theta) \right] & \text{if } \theta \leq t < 3\theta, \\[2mm] Q_\alpha (2Q_\alpha - I)^{k-1} \left[ \dfrac{1}{\rho c} p(t - 2k\theta + \theta) \right. \\[2mm] \left. + (2Q_\alpha - I) \left( \dfrac{2}{\rho c} p(t - 2k\theta - \theta) \right) \right] \\[2mm] \quad \text{if } (2k-1)\theta \leq t < (2k+1)\theta, \ k = 2, 3, \cdots. \end{cases} \tag{43}$$

### 4.2.3. Case 3: $t_1 < 6\theta$

So, $p(t) = 0$ if $t \notin [0, 6\theta)$, and the corresponding exact solutions have the following closed form:

$$v(h,t) = \begin{cases} \dfrac{1}{\rho c} p(t) & \text{if } 0 \le t < 2\theta, \\[2mm] \dfrac{1}{\rho c} p(t) + (2Q_\alpha - I)\left[\dfrac{2}{\rho c} p(t - 2\theta)\right] & \text{if } 2\theta \le t < 4\theta, \\[2mm] (2Q_\alpha - I)^{k-1}\left[\dfrac{1}{\rho c} p(t - 2(k-1)\theta) + (2Q_\alpha - I)\left[\dfrac{2}{\rho c} p(t - 2k\theta)\right.\right. \\[2mm] \left.\left. + (2Q_\alpha - I)\left(\dfrac{2}{\rho c} p(t - 2(k+1)\theta)\right)\right]\right] \\[2mm] \quad \text{if } 2(k+1)\theta \le t < 2(k+2)\theta, k = 1, 2, \cdots; \end{cases} \quad (44)$$

$$v(0,t) = \begin{cases} 0 & \text{if } 0 \le t < \theta, \\[2mm] Q_\alpha\left[\dfrac{2}{\rho c} p(t - \theta)\right] & \text{if } \theta \le t < 3\theta, \\[2mm] Q_\alpha\left[\dfrac{1}{\rho c} p(t - \theta) + (2Q_\alpha - I)\left(\dfrac{2}{\rho c} p(t - 3\theta)\right)\right] & \text{if } 3\theta \le t < 5\theta, \\[2mm] Q_\alpha (2Q_\alpha - I)^{k-1}\left[\dfrac{1}{\rho c} p(t - 2k\theta + \theta)\right. \\[2mm] \left. + (2Q_\alpha - I)\left[\dfrac{2}{\rho c} p(t - 2k\theta - \theta)\right.\right. \\[2mm] \left.\left. + (2Q_\alpha - I)\left(\dfrac{2}{\rho c} p(t - 2k\theta - 3\theta)\right)\right]\right] \\[2mm] \quad \text{if } (2k+3)\theta \le t < (2k+5)\theta, \ k = 1, 2, \cdots. \end{cases} \quad (45)$$

Similar explicit formulas for $t_1 \ge 6\theta$ are excessively cumbersome. In this case, it is easier to directly use the recursive procedure based on (33) and (38) which has the same simple form regardless of the transient load duration and also provides exact results.

## 5. Examples and Discussions

Consider some examples of using the results of the previous section for mathematical modeling of piezoelectric cylindrical devices installed in a car as proposed in [5]. These devices transform the mechanical energy of the moving pistons or crank-shafts into electrical energy, which will be stored in the capacitor or the battery charger. We consider the uniaxial stress state for a cylinder and assume that the material of the cylinder is PZT-5A [4]. In this case, parameters $C, e, \epsilon$ and $\rho$ in Equations (3) have the following values:

$$C = 5.32 * 10^{10} \ \text{N/m}^2; e = 19.89 \ \text{N}/(\text{V} \cdot \text{m});$$
$$\epsilon = 76.12 * 10^{-10} \ \text{farad/m}; \rho = 7750.0 \ \text{kg/m}^3$$

Then the elastic wave speed $c$ in the piezoelectric material is equal to 3684.06 m/s. Next, we take into account that the total force instantaneously applied to the top of a piston in an internal combustion engine is around 6300 pounds, which corresponds to approximately 28,640 N [19]. Suppose that this force $F_a$ is applied downward to a piezoelectric cylinder with a length of $h = 1$ cm and a diameter $d = 1$ cm. So, the area of the upper end face of the cylinder $A = \pi d^2/4 = 0.785$ cm$^2$. Assuming that $F_a$ is uniformly distributed over the upper end face, we get the amplitude of the pressure impulse acting on the top of the cylinder: $p_a = 364.66$ MPA. Let us assume that the applied normal stress load takes the form of the following rectangular compressive load (pressure) impulse:

$$p(t) = -p_a \left[ H(t) - H(t - t_1) \right] \tag{46}$$

where $t_1$ is the duration of the pressure impulse.

We assume first that $\alpha = 0.5$ and $k_\alpha = 1000$ N·s$^{1/2}$·m$^{-1/2}$ (see, e.g., [20]) in the damper boundary conditions (23). Consider the following three values of $t_1$: 5, 10, 15 μs. Since the transit time of elastic waves between the upper and lower end faces $\theta = 2.71$ μs, then these three durations correspond to the three cases of explicit exact solutions considered in 4. The operator $Q_\alpha$ is calculated according to (37). The calculated results for the output voltage $V = \Delta\phi$ (in kV) based on these exact solutions are presented for $t \leq 50$ μs in **Figures 1-3** below.

Comparing the graphs we can see that the maximum or peak value does not depend on the pressure impulse duration $t_1$. However, the number of peaks in each figure depends on the $t_1$. The time distance between two neighboring peaks is approximately equal to $2\theta$. After the pressure load is removed, there is an attenuation of the output voltage vibrations.

Now let us consider another set of the damper parameters: $\alpha = 2.0$ and $k_\alpha = 250$ N·s$^2$·m$^{-2}$. The parameters of the material and the impulse durations are the same as above. The operator $Q_\alpha$ is also calculated according to (37).

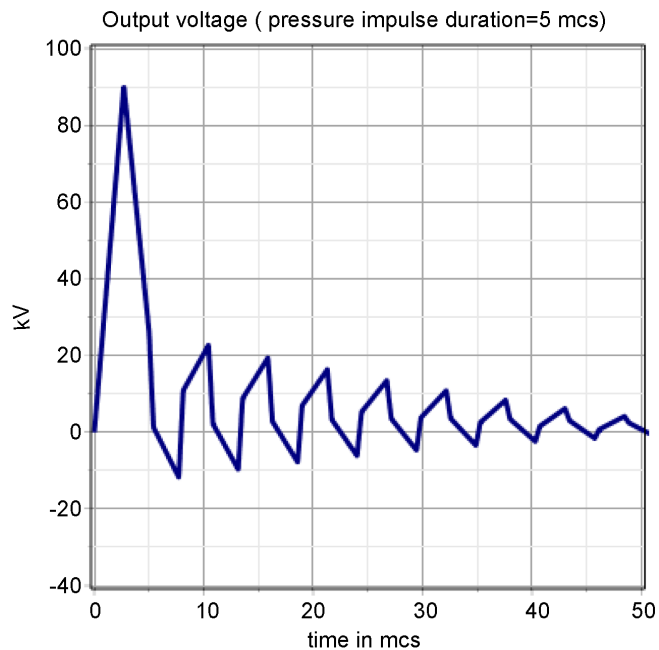

**Figure 1.** Output voltage for $\alpha = 0.5, k_\alpha = 1000$ N·s$^{1/2}$·m$^{-1/2}$ and $t_1 = 5$ μs.

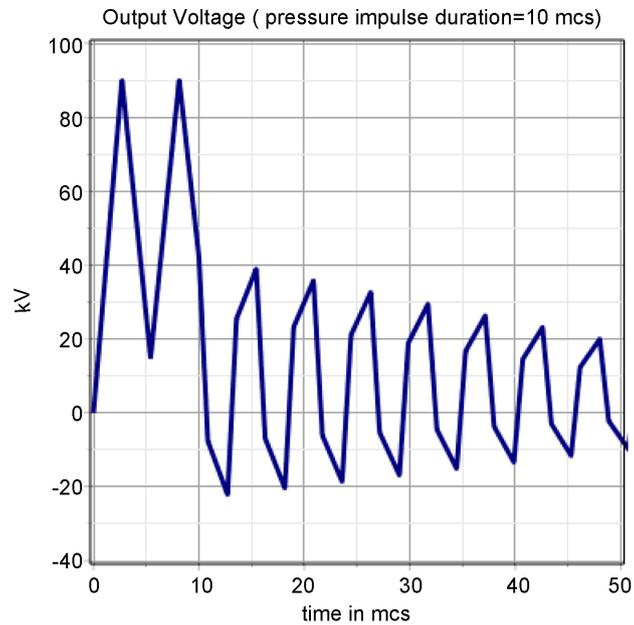

**Figure 2.** Output voltage for $\alpha = 0.5, k_\alpha = 1000 \ \mathrm{N} \cdot \mathrm{s}^{1/2} \cdot \mathrm{m}^{-1/2}$ and $t_1 = 10 \ \mu\mathrm{s}$.

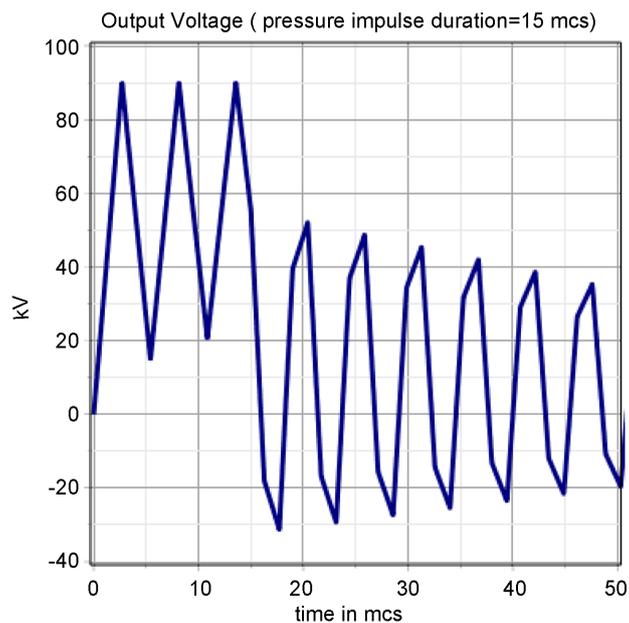

**Figure 3.** Output voltage for $\alpha = 0.5, k_\alpha = 1000 \ \mathrm{N} \cdot \mathrm{s}^{1/2} \cdot \mathrm{m}^{-1/2}$ and $t_1 = 15 \ \mu\mathrm{s}$.

The calculated results for the output voltage $V = \Delta\phi$ (in kV) are presented for $t \leq 50 \ \mu\mathrm{s}$ in **Figures 4-6.**

Comparison of these graphs shows that the maximum value of the output voltage does not depend on the pressure impulse duration $t_1$ which similar to the case when $\alpha = 0.5$. The difference is that now there is only one peak but its width depends on the $t_1$. After the pressure load is removed, the attenuation of the output voltage vibrations is very pronounced: the amplitude of vibrations after the load removal is almost negligible.

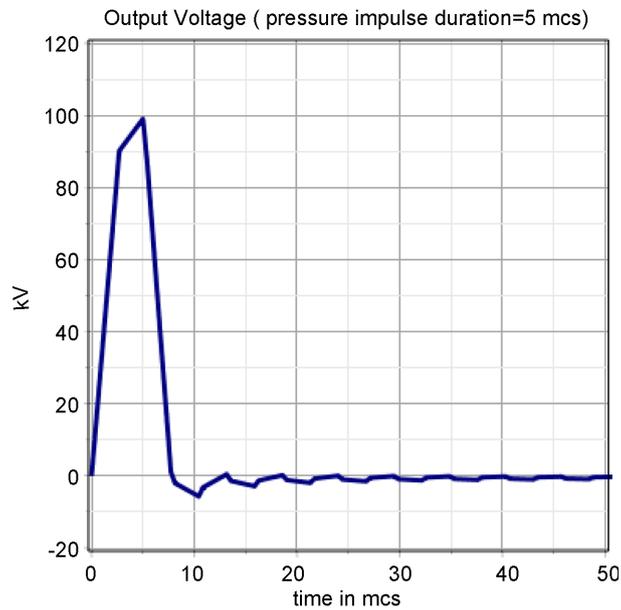

**Figure 4.** Output voltage for $\alpha = 2, k_\alpha = 250 \text{ N} \cdot \text{s}^2 \cdot \text{m}^{-2}$ and $t_1 = 5 \text{ μs}$.

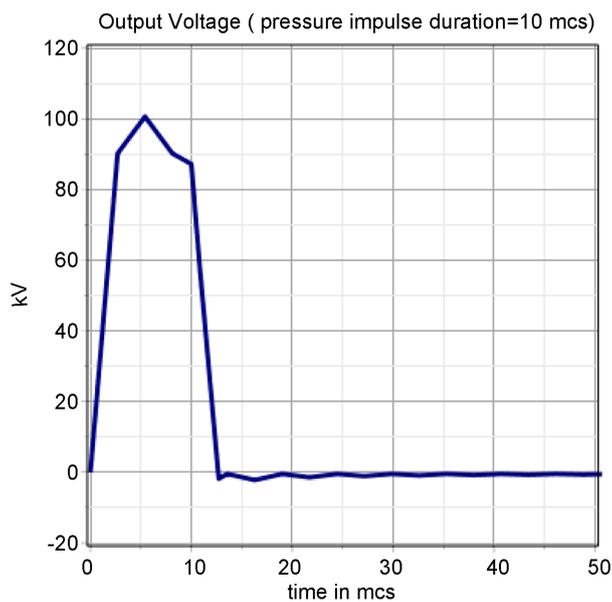

**Figure 5.** Output voltage for $\alpha = 2, k_\alpha = 250 \text{ Ns}^2\text{m}^{-2}$ and $t_1 = 10 \text{ μs}$.

## 6. Conclusion

One-dimensional transient dynamic piezoelectric problems for thickness polarized layers and disks, or length polarized rods, are considered here in the framework of a time-domain Green's function method. As the result, a novel exact analytical recursive procedure is derived which is applicable for a wide variety of boundary conditions including the nonlinear damper case. Some new practically important explicit exact solutions are presented. The effectiveness of the proposed exact approach is demonstrated by examples of the time behavior of the output electric potential difference between two electrodes coated at the end faces of a piezoelectric cylinder fixed to a nonlinear damper at one end, and subjected to impulsive loading at the other.

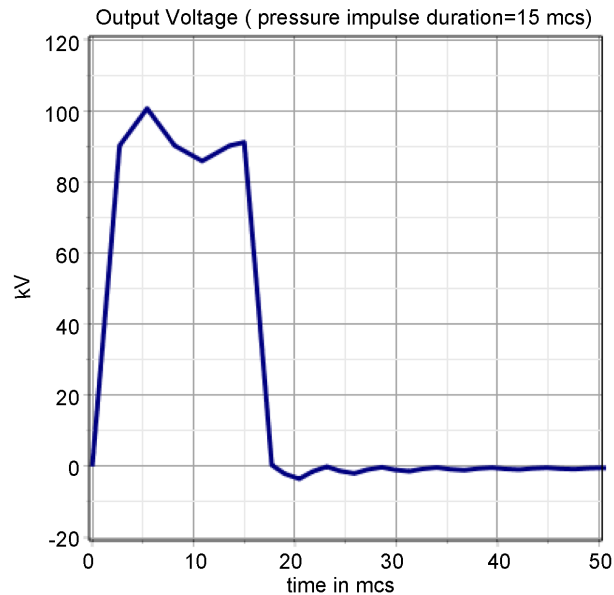

**Figure 6.** Output voltage for $\alpha = 2, k_a = 250 \ \text{N} \cdot \text{s}^2 \cdot \text{m}^{-2}$ and $t_1 = 15 \ \mu\text{s}$.

## Acknowledgements

The authors would like to thank the anonymous reviewers for their comments.